\documentclass[12pt]{amsart}
\usepackage{amsmath,amsthm,epsfig,amssymb}
\usepackage{amsthm,amsfonts,amssymb,amscd}
\usepackage[all]{xy}

\newtheorem{thm}[subsection]{Theorem}

\newtheorem{prop}[subsection]{Proposition}

\newtheorem{lem}[subsection]{Lemma}
\newtheorem{corol}[subsection]{Corollary}
\newtheorem{rem}[subsection]{Remark}

\theoremstyle{definition}
\newtheorem{Def}[subsection]{Definition}

\newtheorem{proposition-definition}[subsection]{Proposition-Definition}

\newcommand{\codim}{\operatorname{codim}}

\newcommand{\ZZ}{{\mathbb Z}}

\newcommand{\PP}{{\mathbb P}}
\newcommand{\NN}{{\mathbb N}}

\renewcommand\square{\frame{\phantom{{\large x}}}}

\author{F. Laytimi}

\address{F. L.: Math\'ematiques - b\^{a}t. M2, Universit\'e Lille 1,
F-59655 Villeneuve d'Ascq Cedex, France}
\email{laytimi@agat.univ-lille1.fr}

\author{W. Nahm}

\address{W. N.: Physikalisches Institut, Universit\"at Bonn,
Nu\ss allee 12, D-53115 Bonn}
\email{werner@th.physik.uni-bonn.de}


\title{A generalization of Le Potier's vanishing theorem}

\begin{document}

\maketitle

\section{Introduction} \setcounter{page}{1}

Consider a vector bundle $E$ of rank $d$  over a compact complex manifold
manifold  $X$ of dimension $n$, and a partition $R = (r_1,r_2,\dots ,r_m)$ of weight 
$r =  \sum\limits_{i=1}^m r_i = |R|, $
where the $r_i$ are strictly positive integers with 
$r_i\geq r_{i+1} $ and $m\leq d$. We call $m$ the length of $R$. 
 
Let ${\wedge} _R = {\mathcal S}_{\tilde R} $
be the Schur functor (for the definition see [M,  p.45]) 
corresponding to the transpose 
$ \tilde R$ of $R$. 
Schur functors were initially  defined on the 
category of vector
spaces and linear maps, but by functoriality 
the definition carries
over to vector bundles on $X$.

\smallskip

We prove the following vanishing theorems for cohomology groups:

\begin{thm} 
 For any partition $R=(r_1,r_2,\dots ,r_m)$, {\ }  
$$H^{p,q}(X, {\wedge}_R E) = 0
\mbox{\quad if\quad  }{\wedge}_R E\ \mbox{is ample\quad  and} 
\quad  p+q-n > 
 \sum\limits_{i=1}^m r_i(d-r_i).$$
\end{thm}

If $m=1$ we get Le Potier vanishing theorem [LP1](see section 2).

\begin{corol}
 
$H^{p,q}(X, \otimes_{i=1}^l\mathcal S^{k_i}E\otimes _{j=1}^m 
\wedge^{s_j} E)= 0$

 if  $\otimes_{i=1}^l\mathcal S^{k_i}E
\otimes _{j=1}^m \wedge^{s_j} E$ is ample and

$p+q-n > {\sum\limits}_{j=1}^m s_j(d-s_j) +(d-1)\sum\limits_{i=1}^l k_i.$

\end{corol}

In [LP2]  Le Potier introduced a very useful tool
 for the derivation of 
such  theorems. 
It is based on the Borel-Le Potier  spectral  sequence, 
a term introduced 
by Demailly [D1].

Given a sequence of integers  $0=s_0<s_1<s_2< \ldots <s_l\leq d,$
and a complex vector space $V$ of dimension $d$,
let ${\mathcal F}l_{s_1,\ldots,s_l}(V) = {\mathcal F}l_{s} (V )$
 be the
variety of partial flags  
$$V_{s_l}\subset V_{s_{l-1}}\subset \ldots \subset V_{s_{1}}
\subset V, \quad 
\codim V_{s_{i}} = s_i.$$

This manifold carries
 canonical vector 
bundles $Q_i$ with fibers $V_{s_{i-1}}/ V_{s_i}$.
 Let 
$Y={\mathcal F}l_{s}(E)$ be the natural fibered variety
 with projection
 
$\pi:Y\rightarrow X$  and fibers ${\mathcal F}l_{s}(E_x),{\ } x \in X.$

We also denote by   $Q_i$  the corresponding vector bundle over $Y$.
 For  partitions ${{\Lambda}_i}$, a vector
 bundle of the form
$\otimes_i\ {\mathcal S}_{{\Lambda}_i}(Q_i)$ over $Y$  will be 
called of Schur type.

The projection  $\pi$  yields a filtration of the bundle 
$\Omega^P_Y$
of exterior differential forms of degree $P$ on $Y$, namely

$$F^p(\Omega^P_Y)=\pi^*\Omega^p_X   \wedge  \Omega^{P-p}_Y.$$  
The corresponding  graded  bundle  is  given  by   
$$F^p(\Omega^P_Y)/F^{p+1}(\Omega^P_Y)=
\pi^*\Omega^p_X\otimes \Omega^{P-p}_{Y/X},$$ 
where $\Omega^{P-p}_{Y/X}$
is the bundle of relative differential forms of degree $P-p$.
For a given line bundle ${\mathcal L}$ over $Y$,
the filtration on $\Omega^P_Y$ induces a filtration on 
 $\Omega^P_Y\otimes {\mathcal L}.$ 
 This latter filtration yields   the Borel-Le Potier
spectral sequence, which abuts to $H^{P,q}(Y,{\mathcal L})$.

It is given by the data $X,Y,{\mathcal L},P$ and will be
 denoted by
$^P{\mathcal E}_B$. Its ${\mathcal E}_1$-terms

$$^P{\mathcal E}_{1,B}^{p,q-p}=  H^q(Y,\pi^*(\Omega^p_X)\otimes
\Omega^{P-p}_{Y/X}\otimes {\mathcal L})$$
can be calculated as limit groups of the Leray
spectral sequence
$^{p,P}{\mathcal E}_L$ associated to the projection $\pi$, 
for which

$$^{p,P}{\mathcal E}_{2,L}^{q-j,j}= H^{p,q-j}(X,R^j \pi _{\ast }{\ }
(\Omega^{P-p}_{Y/X}\otimes {\mathcal L})).$$

For a suitably chosen ample line bundle ${\mathcal L}$ (see section 6) 
of Schur type and for  $P-p = 0$ 
one obtains 

$$^P{\mathcal E}_{1,B}^{P,q-P}=  H^{P,q}(X,{\wedge}_R E).$$

Moreover, under the condition 
$$ (*)  {\ } P+q>n+{\Sigma}  _{i=1}^m r_i(d-r_i),$$ 
the corresponding Borel-le Potier 
spectral sequence will be shown to degenerate at
$^P{\mathcal E}_{1,B}^{P,q-P}$, in the sense that all
 $^Pd_{i,B}$ mapping 
to and from $^P{\mathcal E}_{i,B}^{P,q-P}$ are zero, 
such that
$H^{P,q}(X,\Lambda_R E)$ is a subquotient of 
$H^{P,q}(Y,{\mathcal L}).$
The latter group vanishes by the Kodaira-Akizuki-Nakano
vanishing theorem.

The map $^Pd_{i,B}$  from 
$^P{\mathcal E}_{i,B}^{P,q-P}$ 
is zero by construction, since $F^p(\Omega^P_Y)=0$
for $p>P$. For the maps to $^P{\mathcal E}_{i,B}^{P,q-P}$
we shall use an induction argument to show that their
 sources vanish under the
condition $(*)$.

Because of the Leray spectral sequence, these sources 
are given by subquotients of groups of type $\oplus_j
H^{p,q-j}(X,\ R^j \pi_*(\Omega^{P-p}_{Y/X}\otimes {\mathcal L})).$

Since $\Omega^{P-p}_{Y/X}\otimes {\mathcal L}$ is of
 Schur type, the
vector bundle 
$R^j \pi _{\ast }(\Omega^{P-p}_{Y/X}\otimes {\mathcal L})$
on $X$ is given by a Schur functor applied 
to $E$. This Schur
functor can be calculated for the case where 
$X$ is a point and
$E$ a vector space $V$. Thus the core of our 
proof will be the
evaluation of the groups 
$H^{P-p,j}({\mathcal F}l_{s} (V ),{\mathcal L})$.
 For that, we need a number of
 technical preparations.
\smallskip

In section 3 we explain some basic tools, in 
particular
concerning Young diagrams, which label the
 Schur functors.
In section 4 we reformulate the 
Littlewood-Richardson rules
for the tensor product of two Schur functors. 
This goes slightly 
beyond what we need for the proof, but should 
have independent
interest.
In section 5 we consider the relevant 
cohomology groups on
partial flag varieties, in particular 
on the  Grassmannians.
Section 6 contains the proof of the main
 theorem.
\bigskip

\section{Some known results}

Vanishing theorems for ample vector bundles play a crucial role 
in algebraic 
geometry and its applications.

Let us recall the most important ones.

\begin{thm}[Kodaira - Akizuki - Nakano ] $[KAN]$

Let $L$ be an ample line bundle on an $n$-dimensional
projective manifold $X$.
 
Then $\ H^{p,q}(X,L) = 0 \ \  \mbox{for } p+q > n.$
\end{thm}

The special case $p=n$ is due to Kodaira. In this case
 the ampleness condition
can be relaxed to yield the Kawamata-Viehweg vanishing 
theorem [KV].

The extension of these results to vector bundles 
of higher rank is due  essentially to  Le Potier.

In the sequel,{\ } $E$ is a  vector bundle of rank $d$ on 
a $n$-dimensional
compact complex manifold  $X$.

\begin{thm}[Le Potier] $[LP2]$

If  $E$ is ample, then $H^{n,q}(X,{\wedge} ^r E) = 0 \ \
  \mbox{for } q > d-r.$

\end{thm}

\begin{thm}[Le Potier] $[LP1]$

If ${\wedge} ^r E$ is ample,  then  $H^{p,q}(X,{\wedge} ^r E) = 0
   \ \   \mbox{for } p+q-n > r(d-r).$
\end{thm}

Although Le Potier states  this theorem  under the
 hypothesis that 
$E$ is ample, his proof requires the weaker hypothesis
 that
${{\wedge}}^ r E$ be ample.
Indeed, in his proof he needs $\det\,Q$ on $G_r(E)$ to be 
ample but
$\det\,Q = {\mathcal O}_{\PP ({\wedge}^rE)}(1)|_{G_r(E)}$. 
Thus
 $\wedge^r E$ ample implies that $\det\, Q$ is ample.

\begin{thm}[Sommese] $[S]$

Let $E_j, 1\leq j \leq m$ be ample vector bundles of rank $d_j.$ Then 
$$H^{p,q}(X, \otimes _{j=1}^m \wedge^{s_j} E_j)= 0 \mbox{\quad if\quad  }
\quad  p+q-n > {\sum\limits}_{j=1}^m s_j(d_j-s_j).$$
\end{thm}

Note that Corollary  1.2 gives  Sommese's vanishing theorem 
with $E_j=E$ for all $j$,  under weaker assumption.

For vectors bundles tensored with powers of 
$\det E,$ one obtains less restrictive vanishing conditions. 
The best known case is

\begin{thm}[Griffiths] $[G]$

Let  $E$ be  ample. Then $H^{n,q}(X,\it S ^r E\otimes \det E)) =
 0 \ \  \mbox{for } q > 0.$
\end{thm}

The finite dimensional irreducible  representations of $ Gl(V), $  where
$V$ is vector space of dimension $d$ are in correspondence with partitions of 
length at most $d.$  We denote the irreducible $Gl(V)$-module  
corresponding to the partition $R$ by $S_R(V),$ and call it 
the Schur functor of $V.$

In particular $S^k V= {\mathcal S}_{(k)}V,$ and $\wedge^h V=
{\mathcal S}_{\underbrace {(1,1,\ldots ,1)}_{h\  {\rm times} }}V.$

An extension of the last theorem to Schur functors is

\begin{thm}[Demailly] $[D1]$

Let $R=(r_1, r_2, \dots , r_m)$ be any partition of length $m$.
If $E$ is ample, then $H^{n,q}(X,\it {\mathcal S}_R E\otimes (\det E)^{m}) = 0 
{\ }{for } {\ }q > 0.$
\end{thm}

Note that this theorem can be derived from that of
 Griffiths, as follows:
 
Let $ R = (r_1,\ldots ,r_m)$ be a 
partition of length $m$, and  $r = r_1+\ldots +r_m$
 the weight of $ R$. For $V=\underbrace {E \oplus E \oplus\ldots
  \oplus E}_{m\  {\rm times} },$
 we have   

$\it {\mathcal S}_R E\otimes (\det E)^m \subset  \it S^{r_1}E \otimes \it
S^{r_2}E\otimes \ldots \otimes \it S^{r_m}E \otimes (\det E)^m\subset  
\it S^r V\otimes \det V .$ Then we  use Theorem 2.5.

An extention of the last result to the  whole  Dolbeault cohomology    
is due to Manivel [M2].

\smallskip
For $E$ ample and arbitrary partition ${{\lambda}}$  the question 
of finding an exact condition for $H^{p,q}(X,\it {\mathcal S}_{{\lambda}} E)$
 to vanish  
is still open. For a precise conjecture in the case $p=n$ see [L].

In [LN], as a special case of a more general result, this conjecture 
is proved for any $(p,q)$ and any hook Schur functor
${\Gamma}^\alpha_k E.$  

The latter are defined for $ 0\leq \alpha< k$ and correspond to
 the partition
$(\alpha  +1,1,\dots ,1)$ of weight $k\in \NN^*$. 
Inductively, they can
 be defined as follows:  
$$\Gamma^0_k E={\wedge}^kE$$
and
$${\wedge}^{k-\alpha }E\otimes S^{\alpha } E=
\Gamma ^{\alpha-1}_k E\oplus \Gamma ^\alpha_k E$$
for $0<\alpha<k$.
In particular,  $\Gamma^{k-1}_k E=S^kE$.
Note that $\Gamma ^\alpha_k E=0$\\
 for $d-k+\alpha<0.$

Define a function $\delta:\NN\rightarrow  \NN^*$ by: 
$${\delta(x)\choose 2}\le x<{\delta(x)+1\choose 2}\ .$$
In other words,
$$ 
\begin{array}{llllll}
\delta(0)&=1.\\
\delta(1)&=\delta(2)&=2\\
\delta(3)&=\delta(4)&=\delta(5)&=3\\
\delta(6)&=\delta(7)&=\delta(8)&=\delta(9)&=4\\
\end{array}
$$
etc...
\medskip

\begin{thm}  If $E$ is ample,  then
$H^{p,q}(X,\Gamma ^\alpha_k E)=0 $ for

$q+p-n>(\delta(n-p)+\alpha)(d-k+2\alpha)- \alpha(\alpha+1)\ .$
\end{thm}
For $\alpha=k-d$ one obtains
$H^{p,q}(X,S^\alpha E\otimes \det E)=0,$
when
$p+q-n>(\delta (n-p)-1)(k-d)$.
For $p=n$, this specializes to Griffiths' vanishing theorem.

In the present paper we prove an extension to Schur
 functors  of Le Potier's
theorem 
(Theorem 2.3~).

\section{Basics definitions and tools}

\subsection{Some notations and definitions}$\phantom{0}$\\

$\NN^*=\NN-\{0\},{\ } {\ } I(r) = \{1,2,\dots, r\}
 \subset \NN^*.$
\smallskip

For $(i,j) \in \ZZ \times \ZZ$, we call $i$
the height and $j$ the width of $(i,j)$.\\
For $S \subset I(r) \times \mathbb Z, {\ } {\rm card} (S) < 
\infty,$ {\ }
  we define
the sequence{\ } $$[S] : I(r) \longrightarrow \mathbb N {\ } 
{\mbox by }{\ }
[S]_i = {\rm {\rm card}}  \{j \in \mathbb Z {\ }|{\ } (i,j)
 \in S \}.$$
For $S \subset I(r) \times \mathbb N^*, {\ } {\rm {\rm card}} 
(S) < \infty,$
we define $$\langle S \rangle: \mathbb N^* \rightarrow 
\mathbb N
{\ }{\mbox by }{\ }\langle S \rangle_j = {\rm {\rm card}} 
\{i \in  I(r) {\ }|{\ }(i,j)
 \in  S \}.$$
\medskip

Partitions $u$ of lenght  $l(u)\leq r$ are  
weakly decreasing sequences 
 $u: I(r)\rightarrow \mathbb N.$ More precisely, 
 we regard partitions
$u= (u_1,u_2,  \ldots,u_r  )$ {\ } and 
$(u_1,u_2,  \ldots,u_r, 0, \ldots, 0)$ as equivalent. For
$u_r\neq 0$ the length of $u$ is $r$.

The weight $u$ is given by 
$$|u|=  {\rm {\rm card}}( Y(u))$$
where
$$Y(u) = \{(i,j)\in I(r)\times \NN^*{\ } | {\ } 1 \leq j\leq
u_i\}$$
is the Young diagram of $u$. Equivalently,
$|u|= \sum\limits_{i} u_i$.  

For example the Young diagram of $u=(4,2,1,0)$ is 
$$Y(u) = \{(1,1),(1,2),(1,3),(1,4),(2,1),(2,2),(3,1)\},$$
the length is 3 and the weight is $|u|=7.$

The transpose $\tilde u$ of a partition $u$ is defined by

 $$Y(\tilde u) = \widetilde {Y (u)}, $$
where $\widetilde {(i,j)} = (j,i)$ for $(i,j)\in\NN^*\times
 \NN^*.$
For non positive integers $ j$ we put $\tilde u_j=+\infty$.

We define the squared norm of a partition $ u$
by  $||u||^2 = \sum \limits_ {i\in \NN^*}\tilde u_i^2.$
\medskip

When $\psi$ is a map $I(r)\rightarrow \mathbb Z,$ we denote the
corresponding weakly decreasing sequence by $\psi^{\geq}.$
More precisely, 
$$\psi^{\geq}= \psi \circ\sigma,$$
where $\sigma$ is any permutation of $I(r)$ such that
 $\psi\circ\sigma$
is weakly decreasing. When the image of $\psi$ lies in $\mathbb N^*$, 
$\psi^{\geq}$ is a partition of length $r$.
\smallskip

When $u$ is a weakly decreasing sequence, let $u^>$ be
 the sequence
obtained by removing repetitions, in other words the 
strictly decreasing
sequence which has the same set of terms as $u$. For 
a finite sequence
$u$ let $u^<$ be the strictly increasing sequence 
which has the same
set of terms as $u$.

Every partition $u$ can be reconstructed from $a=u^>$ 
and $s=\tilde u^<$.
We write $u=a_s$. Explicitly
$$a_s=(\underbrace {a_1,\ldots ,a_1}_{s_1\ {\rm times}},
\underbrace {a_2, \ldots ,a_2}_{(s_2-s_1)\ {\rm times}}
,\ldots,
 \underbrace {a_j, \ldots ,a_j}_{(s_j-s_{j-1})\ {\rm times}}
 ,\ldots).$$
The latter notation also will be used for general sequences 
$a$ of finite length.

\medskip

For each partition $u$ one has a Schur functor 
${\mathcal S}_u :{\mathcal V}
 \rightarrow \mathcal V,$
where $ \mathcal V$ is the category of complex vector spaces 
and linear maps.
If $u_r>0$ for some $r>\dim{\ }V$, one has 
$\ {\mathcal S}_u V = 0.$

By functoriality ${\mathcal S}_u$ also operates on 
vector bundles
over a given manifold $X.$

A generalized partition of length $r$ is a weakly 
decreasing sequence {\ }$ u:I(r) \rightarrow \mathbb Z.$ {\ }

We define the diagram of a generalized partition 
$u$   by
$${\mathcal D}(u) = \{(i,j)\in I(r)\times \mathbb Z {\ }|{\ } 
j{\ }\leq u_i \}.$$

Note that this diagram is infinite even when $u$ is a 
usual partition.

\smallskip

{\bf {Example}}

\smallskip

For the generalized partition $u = (2,1,0,-1,-2),{\ } {\mathcal D}(u) $ 
is the set of the following marked  points :
  
\smallskip

\begin{figure}[htb]\label{fig2}
\begin{center}
\epsfig{file=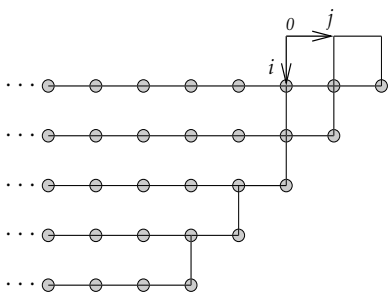,height=1.5in,width=2in,angle=0.0}
\caption{diagram of a partition}
\end{center}
\end{figure}

We define the involution 
$$\chi ^*:I(r)\times \mathbb Z \rightarrow I(r)\times \mathbb Z {\ }{\ }
{\mbox  by} {\ }{\ }\chi ^*(i,j) = (r+1-i,\ 1-j),$$  and
the reversed generalized partition $\chi(u)$ 
  by

$$ {\mathcal D}(\chi (u))^c = \chi^*({\mathcal D}(u))$$
where $({\ } )^c$ denotes the  complement in $I(r)
\times \mathbb Z.$
Explicitly
 
$$\chi(u) = (-u_r, \dots, -u_2,-u_1){\ }{\ } {\mbox {\rm for}} {\ }{\ }
u = ( u_1,u_2, \dots,u_r).$$
\smallskip

For the category  ${\mathcal V}_r$ of complex vector spaces of
 dimension $r$ we
extend the Schur functor notation 
${\mathcal V}_r \rightarrow {\mathcal V}$
to generalized partitions of length  $r$ by
$${\mathcal S}_{u-{\bf k}_r}{ V}  =  {\mathcal S}_u{ V} \otimes 
(\det{\ } {V}^*)^k,
\quad k\in\NN^*,$$
where ${\bf k}_r$ is the partition $(k,k,\dots, k)$ 
of length $r$
and $V$ is a complex vector space of dimension $r.$

For $u = ( u_1,u_2, \dots,u_r)$ we have 
$${\mathcal S}_u{ V}^*\simeq {\mathcal S}_{\chi(u)}V.$$

If a sequence $u$ is not weakly decreasing, we put 
${\mathcal S}_uV=0$.
\bigskip

If {\ }  ${\mathcal D}(v)\subset {\mathcal D}(u)$,{\ }
   we say that the 
pair $u,v$ forms
a skew partition $u/v$.
We define the diagram of such a skew partition by
$${\mathcal D}(u/v) = {\mathcal D}(u){\ } /{\ }
 {\mathcal D}(v)$$

where ${\mathcal D}(u)\,/ \,{\mathcal D}(v)$ denotes 
difference of sets,
 
 and the weight by
$$|u/v| = {\rm {\rm card}} ({\mathcal D}(u/v))$$ 
 
The reversed skew partition is given by
$$\chi (u/v) = {\chi}(v) / {\chi}(u)$$

\subsection{On the dominance partial order and ampleness}

\begin{Def}
Let $I=(i_1,i_2\dots )$, $J=(j_1,j_2\dots )$ be partitions 
 of the same weight. We say that
$$I\succeq J {\ } {\ }  {\ }if {\ } {\ } {\ }\sum
 \limits _{k=1}^l 
{i}_k{\ }\geq {\ }\sum \limits_{k=1} ^l {j}_k {\ }  {\ }{\mbox
{for {\ } any}}{\ }{\ } l.$$
This relation is called the dominance partial order.
\end{Def}

\smallskip

\begin{Def}
For partitions $I$, $J$ of arbitrary weight, this definition  
is generalized in [LN]
 by $I {\ }\preceq {\ } J$ for
$|J|I{\ }\preceq {\ } |I|J$. Here the multiplication
 of a partition $I$
by $n\in\NN$ is defined by $n(i_1,i_2,\ldots)=(ni_1,ni_2,\ldots)$.
\end{Def}

\smallskip

{\bf {Example}}

\smallskip

For the partition of weight 5{\ } $I = (1,1,1,1,1), $
and {\ }   $J = (2,1)$ of weight 3, we have $I \preceq  J$ because
 $1/5 < 2/3,{\ }
2/5 < 1, {\ } 3/5 < 1$ and $4/5 < 1.$

The following lemma of Macdonald concern partitions of the same 
weight, but his proof can be adapted easily to partitions of arbitrary weight.

\begin{lem}$[M, p.7]$ 

\smallskip

For any non-trivial partitions of arbitrary weight $I$,  $J$

$$I {\ }\preceq {\ } J \Longleftrightarrow \tilde I {\ } \succeq {\ } 
\tilde J.$$

\end{lem}

In [LN] we also proved   
 
\smallskip

 \begin{lem}$[LN]$

\smallskip
 
For any partitions $I$ and $J:$  

$$\mbox{if} {\ } {\ } I \succeq J,  {\ }\mbox {then}{\ } {\ }
{\mathcal S}_IE{\ } \mbox {ample}
 \Longrightarrow {\mathcal S}_JE {\ }  \mbox {ample}.$$

\smallskip

In particular,

$$\mbox{if} {\ } {\ }I \simeq J {\ } {\ }{\ }  \mbox 
{then} {\ } {\ }{\ }
{\mathcal S}_IE {\ }{\ } \mbox {ample}\Longleftrightarrow
{\mathcal S}_JE {\ }{\ } \mbox {ample}.$$
\end{lem}
 
 We write $I \simeq J $ if $I \succeq J $ and $I \preceq J.$ 
For example $(k,0,0,\ldots )\simeq (1,0,0,\ldots)$

\section{On the Littlewood-Richardson rules}

\begin{Def}
On any skew  partition $w/u$, we define on   ${\mathcal D} (w/u)$ 
the Littlewood-Richardson 
(LR) order by

$$(i, j)<_{LR}(i', j') {\ } {\mbox{for {\ }}}  {\ } i<i'$$
$$(i, j)<_{LR}(i, j') {\ } {\mbox {for {\ }}} {\ } j>j'.$$

In this section $w,u$ are fixed and $x<_{LR}y$ implies 
$x,y \in {\mathcal D} (w/u)$.

Let $w/u$ be a skew partition, and  $v = (v_1,\dots, v_r)$
a partition such that $|v|=|w/u|.$ We denote $c_1$ the numbering

$$c_1  : {\mathcal D}(w/u)\rightarrow \NN^*
{\ } \mbox {with}{\ } 
{\rm {\rm card}} \{x\in  {\mathcal D} (w/u){\ } | {\ } c_1(x)=
k\}=v_k,\  \forall k\in \NN^*.$$
\end{Def}

\begin{Def}\label{LR rules}

The  numbering $c_1$ 
is said to satisfy the LR rules , iff

$(L_1)$ : $c_1$ is strictly increasing on each column of
 ${\mathcal D} (w/u),$

$(L_2)$ : $c_1$ is weakly increasing on each row of  
${\mathcal D} (w/u),$

$(L_3)$ : for all $x\in {\mathcal D} (w/u)$ and  all 
$k\in \NN^*$, one has
$\sigma_k(x)\geq \sigma_{k+1}(x), $ where 
 
$\sigma_k(x)={\rm {\rm card}} \{y\leq_{LR} x {\ }|{\ } c_1(y)=k \}$

\end{Def}

\begin{Def}
For a given $c_1$, let $c_2 : {\mathcal D} (w/u)\rightarrow 
\NN^*$ {\ }
with {\ } $c_2(x) = {\rm {\rm card}} \{y\leq_{LR} 
x\ |\ c_1(y) = c_1(x)\}.$

Then $c=(c_1,c_2),{\ } c:{\mathcal D} (w/u)\stackrel{\sim}
{\longrightarrow} Y(v)$ is  a bijection. 
Let $b=c^{-1},$ with $b=(b_1, b_2).$

\smallskip

\begin{lem}

The map $b$ is a strictly increasing
function of the width  on each row of $Y(v)$
with respect to the order $<_{LR}.$

\end{lem}

{\it Proof:}

Let $x=(i,j)$, $y=(i,j')$, $j>j'$. Then
${\rm card}\{z\leq_{LR}x\ |\ c_1(z) = i\} = j > j' =
{\rm card}\{z\leq_{LR}y\ |\ c_1(z) = i\}$,
thus $b(i,j)>_{LR} b(i,j')$.
$\hfill{\square}$

We put  $C(x)=\{c(y)  {\ }|{\ }  y \leq_{LR} x \}$.

\end{Def} 

\begin{rem}

$(i,j)\in C(x) \Longleftrightarrow b(i,j)\leq_{LR} x.$  
By the previous lemma,
$i,j \in \NN^*, (i,j)\in C(x)  \Longleftrightarrow j\leq 
{\rm {\rm card}} \{y\leq_{LR} x{\ }|{\ }  c_1(y) = i \}
=\sigma_i(x)$.
\end{rem}
\smallskip

\begin{lem}
   
Assume that $c_1 $ satisfies the LR conditions
for each $x \in {\mathcal D} (w/u)$. Then  all $C(x) $
 are Young diagrams. 
\end{lem}
 
{\it Proof:} 

We  have to show $(i,j)\in C(x) \Longrightarrow (i',j) \in C(x) {\ } 
\mbox {for} {\ } i'=1, \dots, i-1,$ and 
$(i,j)\in C(x) \Longrightarrow (i,j') \in C(x) {\ } \mbox {for}
 {\ }j'=1, \dots, j-1$. 
The first implication follows from
 
\smallskip

$j\leq \sigma_i(x) \leq^{(L_3)} \sigma_{i'}(x)$.
The superscript over the inequality
sign  indicates that this inequality holds by
 virtue of the corresponding property.
  
The second implication follows from remark 4.5.  
$\hfill{\square}$

\smallskip
\begin{Def}

We say that the set $C(x)$ satisfies the condition  $(Y)$,  

If for each $x \in {\mathcal D} (w/u),\ C(x) $ is a Young diagram.

\end{Def}

\smallskip

\begin{lem}
\smallskip
Let $c = (c_1,c'_2),{\ } c :{\mathcal D} (w/u)\longrightarrow Y(v)$
  be 
 a bijection, such that $(Y)$ is true and $c_1$ satisfies 
 $(L_1),(L_2).$ 
Then $c_1$ satisfies the LR rules and $c'_2 = c_2.$
\end{lem}

{\it Proof:}

The first part of the conclusion is obvious,
since $C(x)$ is a Young diagram, and the length of the
$k$-th row of $C(x)$ is equal to  $\sigma_k(x).$

The restriction of $c'_2$ to the 
set
$$\{x\in {\mathcal D} (w/u) {\ }|
\break {\ } c_1(x)=i\}$$ is an order-preserving
map to the $i$-th row of $Y(v)$, where the latter
is ordered by the width, Indeed 

For $x,y \in {\mathcal D} (w/u)$ with $c_1(x) = c_1(y)=i,$ we have 
$x = b(c(x))= b(i,c'_2(x)),$ and $y = b(i,c'_2(y)), $ then 
$ c'_2(x) < c'_2(y) \Longleftrightarrow x<_{LR}y$ by the Lemma 4.7. 

Since both sets have $v_i$ elements, this map is the unique
order-preserving bijection.

$\hfill{\square}$

\smallskip

From now on we will use the

\begin{Def}

We say that a bijection 
$c:{\mathcal D} (w/u) \stackrel{\sim} {\longrightarrow} Y(v)$ 
satisfies the LR rules, iff $(L_1), (L_2)$ of the definition 
4.2 and (Y) are satisfied.

\end{Def}

 Recall that $b=c^{-1}$.
\smallskip

The importance of the LR rules is due to the
following well-known proposition.

\begin{prop}

Let $\dim{\ } V = r, {\ }{\ }u $ a generalized partition of
 length $r$ and $u'$ a partition. One has 

$${\mathcal S}_u V \otimes {\mathcal S}_{u'} V \simeq \bigoplus_
{(w,b)\in LR(u,u')} {\mathcal S}_w V,$$
where $LR(u,u')$ consists of pairs $(w,b)$ such that $w/u$ is a skew
partition and 
$$b: Y(u')\stackrel {\sim}\longrightarrow {\mathcal D}(w/u)$$
a bijection satisfying the Littlewood-Richardson rules.
\end{prop}

\smallskip

The LR rules have a useful symmetry which is hidden in their
original definition, but will be made explicit in the following
proposition.

\begin{prop} 

\smallskip

$c$ satisfies the LR rules , iff

\begin{description}

\item [(h)]   On each column of ${\mathcal D} (w/u), c_1$ 
preserves the
 order of the heights 

\item [(th)]    On each column of ${\mathcal D} (w/u),c_2 $ weakly 
inverts this order

\item [(w)]    On each row of ${\mathcal D} (w/u),  c_2$  inverts
 the  order 
of the widths 

\item [(tw)]    On each row of ${\mathcal D} (w/u), c_1$  weakly 
preserves
 this order 

\item [(h' )]    On each column of $Y(v), b_1 $ preserves 
the order of the 
 heights 

\item [(th')]    On each column of $Y(v), b_2$  weakly 
inverts this 
order 

\item [(w')]    On each row of $Y(v), b_2$  inverts the  order
 of the widths

\item [(tw')]    On each row of $Y(v), b_1$  weakly 
preserves this order.
 
\end{description}

\end{prop}

{\it Proof:}

Conditions $(\rm h), (\rm tw)$ are restatements of $(L_1),(L_2)$,
and $(\rm tw')$ follows from lemma 4.4.

To show $(L_1,L_2, Y) \Longrightarrow (\rm th)$ we show 
$(L_1, L_2, L_3)
\Longrightarrow (\rm th),$ 
which is equivalent by the lemmas 4.6 and 4.8.

Let $x,y,x',y'\in {\mathcal D} (w/u)$ with $y = (i-1, j),{\ }
 x = (i, j)$  

and $y' = (i-1, j-1),{\ }  x' = (i, j-1), $ with $c_1(y) = k,
{\ }     
c_1(x) = k',$ by $(L_1), k < k'$. We have $\sigma_k (x)
\leq \sigma_k (y). $

Consider  the case

a){\ } $\sigma_k(x) =  \sigma_k(y) $

This case corresponds to $y'\notin {\mathcal D} (w/u)$ or 
$c_1(y')< c_1(x)-1.$ We
have by $(L_3){\ }  \sigma_{k'}(x) \leq \sigma_{k}(x) $ 
hence $\sigma_{k'}(x) \leq
\sigma_{k}(y), $ which is by definition $c_2(x)\leq c_2(y).$ 

\smallskip
b){\ } $\sigma_k(x) > \sigma_k(y) $

This case corresponds to $y'\in {\mathcal D} (w/u)$ and 
$c_1(y')\geq c_1(x)-1.$ 

Moreover

$c_1(y')<^{(h)}c_1(x')\leq^{(tw)}c_1(x)${\ } and {\ }  
 $c_1(y')\leq^{(tw)}c_1(y)<^{(h)}c_1(x).$

This implies 

$c_1(x) = c_1(x')$ and $ c_1(y') = c_1(y) = c_1(x) -1.$

This gives $\sigma_{k'}(x) = \sigma_{k'}(x')-1$ and 
$\sigma_k(y) = \sigma_k(y')-1. $   

Now we use the induction on $j$, so that we can assume
 $c_2(x') \leq c_2(y'). $ Thus 

$c_2(x)=\sigma_{k'}(x) = c_2(x')-1$ and $c_2(y)=
\sigma_{k}(y) = c_2(y')-1. $

The starting step of the induction is when $x$ and
 $y$ are such that
  $y'\notin {\mathcal D} (w/u), $ which is the case of a).  

To prove $(\rm tw),(Y) \Longrightarrow (w),$  assume 
$(i, j), {\ }  (i, j+1) \in {\mathcal D} (w/u)$.
We want to show that the assumption 
$c_2(i, j+1)\geq c_2(i, j)$ 
leads to a contradiction.

Now $c_1(i, j+1)\geq^{(tw)} c_1(i, j)$,   
and the two inequalities together imply that $c(i, j)$
 belongs to  
the Young diagram $C((i, j+1)),$ which is wrong, 
since $(i, j) >_{LR} (i, j+1).$

\smallskip

To prove $(\rm h')$, assume that $i'<i, {\ } b_1(i',j)
\geq b_1(i, j).$

The case $b_1(i',j) =  b_1(i, j)$ is excluded by $(\rm w)$,
 which just has been proved.
For $b_1(i', j) >  b_1(i, j),$ we have
 $(i', j)\notin C(b(i, j)), (i, j) \in C(b(i, j)),$
which contradicts $(\rm Y).$

\smallskip

To prove $ (\rm h'), (\rm w), (\rm th) \Longrightarrow 
(\rm th')$ assume for $x= (i,j), x' = (i',j)$, $x,x'\in Y(v)$
 that $i'<i$, $b_2(x) > b_2(x')$. Let $z= (b_1(x), b_2(x'))$.
Since $b_1(x) > b_1(x')$  by $(\rm h')$ we have
$z\in {\mathcal D} (w/u)$.
This yields

\smallskip 
$j = c_2(b(x))\stackrel{(w)}<c_2(z)\stackrel {(th)}
\leq c_2(b(x')) = j,$ which is absurd.

To prove 
$(\rm h),(\rm tw),(\rm Y) 
\Longrightarrow (\rm w')$ let $(i, j'), (i, j)
 \in Y(v),
j'<j.$ By lemma 4.4 we have $b_1(i, j')\leq b_1(i, j)$.

In the case $b_1(i, j') = b_1 (i, j),$  we have 
 $ b_2 (i, j) < b_2(i, j'),$ by lemma 4.4,
thus $(\rm w').$
To prove $(\rm w')$ in the case $b_1(i, j')< b_1(i, j)$,
assume $b_2(i, j') \leq b_2(i, j).$ Then $z \in  {\mathcal D} 
(w/u),$ 
where $z = ( b_1(i, j),b_2(i, j')).$

Thus $i = c_1(b(i, j')) \stackrel {(h)}<c_1(z) 
\stackrel {(tw)}\leq c_1(b(i, j)) = i,$ which is absurd.

\smallskip

Finally let us show
 $(\rm h'), (\rm w'),(\rm tw') \Longrightarrow (\rm Y).$
The implication

$(i, j)\in C(x),{\ } j'<j \Longrightarrow (i, j')\in C(x)$
 is  equivalent  to 
  
 $(j'<j \Longrightarrow b(i, j')<_{LR} b(i, j)).$

For $b_1(i, j') =  b_1(i, j)$  this  follows from $(\rm w')$  
(or from $(\rm w).$) 
The case $ b_1(i, j)<_{LR} b_1(i, j')$ is excluded by $(\rm tw').$

The implication

$(i, j)\in C(x), i'<i \Longrightarrow (i', j) \in C(x)$
 is  equivalent  to
  
$(i'<i \Longrightarrow b(i', j) <_{LR} b(i, j)).$

By $(\rm h')$
one even has {\ } $b_1(i', j) < b_1(i, j).$

$\hfill{\square}$

\begin{rem}
Clearly the conditions $(\rm h, \rm w,\rm th,\rm tw,\rm h'
,\rm w',\rm th',\rm tw')$ 
are not independent.
For example, we have proved
$(\rm h,\rm tw,\rm w',\rm h',\rm tw')\Rightarrow (\rm h,\rm tw,Y)$ 
$\Rightarrow (\rm w,\rm th,\rm th').$
\end{rem}
\smallskip

\begin{rem}
The set of conditions remains invariant under the replacement

$$ Y(v) \rightarrow Y(\tilde v), \qquad \mbox {and} \qquad
{\mathcal D} (w/u)\rightarrow \widetilde{\chi^*({\mathcal D}(w/u))},$$

which exchanges columns and rows. 
\end{rem}

\begin{Def}

For $U\subseteq \ZZ \times \ZZ$, a map

 $b:U  \longrightarrow \ZZ\times \ZZ $ with 
 $b(i,j) = (b_1(i,j), b_2(i,j))$
is called height increasing if 
$b_1(i,j)\geq i,{\ } \forall (i,j)\in U $
\end{Def}
\begin{rem}
For partitions $u,v$ of the same weight,
 the dominance partial
order can be characterized by the property that  
$u\preceq v$ if there is a height increasing bijection 
$b : Y(v) \longrightarrow Y(u).$
\end{rem}
\smallskip

\begin{lem}

If $A\subset \NN^*\times \NN^*$ and 
$b:A \longrightarrow \NN^*\times \NN^*$ 
preserves the order of the height on each
column of $A$ , then $b$ is height increasing.

\end{lem}

{\it Proof:}

By induction, since 
$1\leq b_1(1,j) < b_1(2,j) < \ldots < b_1(i,j) $
for all $(i,j)\in A.$
This implies {\ } $b_1(2,j)\geq 2$ {\ } etc $\ldots$

\section{Cohomology  groups on  flag manifolds}

\medskip

For $V$ a vector space of dimension $d$ and a sequence
$s=(s_0,s_1,\dots , s_l)$ such that
$0 = s_0<s_1<s_2<\ldots <s_l<d$,
the  flag manifold 
${\mathcal F}l_{s} (V ) = {\mathcal F}l_{s_1,s_2,\ldots,s_l}(V)$
given by subspaces  $V_{s_i}\subset V$ of codimension $s_i$
has natural vector bundles $Q_i$ with fibers 
$V_{s_{i-1}}/ V_{s_i}$.

For a partition $a=( a_1,a_2,\ldots,a_l)$ such that 
$a_1>a_2> \ldots >a_l$,
we consider the Schur type line bundle
$$Q^a=\bigotimes_{k=1}^l {\ } \det (Q_k)^{a_k}\ .$$

Our aim is to prove

\begin{thm}
\[  H^{p,q}({\mathcal F}_s(V), Q^a) = \left\{  \begin{array}{llll} 
\delta_{q,0}\ {\mathcal S}_{a_s} V  &\mbox{if } {\ } p = 0 \\ 
\bigoplus_{i\in I(p,q,s,a)}\ {\mathcal S}_{\rho(i)} V 
&{\mbox{if }} {\ } p\neq 0,
\end{array} \right. \]
where $|\rho(i)| = |a_s|$\ and for all $i$ in the index set
$I(p,q,s,a),$
\begin{description}

 \item [(i)]  
$$\rho(i) \prec a_s$$ 

\item [(ii)] 
$$p+q+1+||a_s||^2 \leq ||\rho(i)||^2.$$
\end{description}
\end{thm}

The main tool for deriving such results is 
\bigskip

{\bf Bott's Theorem :} [ B ] or [D2]

\smallskip

{\it Let $V$ be a complex vector space of dimension $d$ and
 ${\mathcal F}(V)$
the complete flag manifold of $V$. 
Let $a \in \ZZ^d$ and $I(d)= (1, 2, \dots ,d)$. {\ }

Define $\psi(a)= (a-I(d))^{\geq} + I(d)$,
 where $(a-I(d))^{\geq}$
is the sequence obtained by rearranging the terms of $(a-I(d))$
in weakly decreasing order.

We call {\ } $i(a)$ the number of strict inversions of $(a-I(d))$: 

$i(a) = {\rm {\rm card}}\{(i,j){\ }|{\ }i<j,\
 (a-I(d))_i < (a-I(d))_j \}.$}

Then

$$  H^q({\mathcal F}(V), Q^a) = \delta_ {q, i(a)}  
{\mathcal S}_{\psi (a)}V .$$

\bigskip

Recall that one puts ${\mathcal S}_{\psi (a)}V=0$ if ${\psi (a)}$
 is not a
partition.
In particular the cohomology of $Q^a$ is non vanishing  iff all
components of {\ } $a-I(d)$ {\ }  are pairwise  distinct.

\begin{Def}For $a \in \ZZ^d$, we say that {\ } $a$ {\ }is
 admissible iff all
components of 
{\ }$a-I(d)$ {\ }  are pairwise  distinct. 
\end{Def}

\begin{corol}

For the incomplete flag   ${\mathcal F}l_s(V)$, we have [M1]

$$H^q({\mathcal F}l_s(V), Q^a) = \delta_ {q, i(a)} 
 {\mathcal S}_{\psi (a_s)}V
.$$
\end{corol}

 \smallskip
 
 Let $G_r(V)$ be the Grasmannian manifold of codimension $r$
 subspaces of $V,$
$Q$ and $S$ the  universal quotient bundle and the
 universal subbundle 
 on 
this manifold. Then the following statement holds:

\begin{corol} 

For $u,v$ generalized partitions of lengths $r$ and $d-r$ respectively, 
we have [M1]

$$H^q(G_r(V),{\mathcal S}_u Q \otimes {\mathcal S}_v S) = \left\{
  \begin{array}{crcr}
 {\mathcal S}_{\psi(u,v)}V &{\mbox if}{\ } q = i(u,v),\\
 0 & {\mbox otherwise}. 
 \end{array} \right. $$

\end{corol}
 
 Note that $|\psi(u,v)|=|u|+|v|$. 
 
Let $w$ be a generalized partition of length $r$ and
 $u$ a partition, let $a = (w, \tilde u ), $ then the set of elements in 
$(a-I(d))$ 
is $\{ \{\alpha_i\}_{i=1, \ldots, r},
\{ \beta_j\}_{j=1, 
\ldots,d-r} \}$, where  ${\ } \alpha_i = w_i-i ,{\ } \beta_j = 
\tilde u_j-(r+j).$  

$(w, \tilde u ) $ is admissible iff 
$ \forall (i,j) \in I(r)\times I(d-r), \alpha_i \neq \beta_j.$
We have $i(a)  = {\rm {\rm card}}\{(i,j){\ }|{\ } \alpha_i <  \beta_j \}$.

Let $[\gamma]$ be the sequence such that, 
$[\gamma]_i = {\rm card}\{j{\ }|{\ }
 \alpha_i-\beta_j < 0 \}$,\break 
and  $\langle \gamma\rangle$ the sequence such that, 
 $\langle \gamma\rangle_j = {\rm card} \{i{\ }|
{\ }\alpha_i-\beta_j < 0 \}$.

\begin{lem}
 
Let $w$ be a generalized partition of length $r$ and
 $u$ a partition,
such that $(w, \tilde u)$ is admissible. 
Define  $s_+ : I(r) \rightarrow \ZZ${\ } by {\ }
$(s_+)_i =w_i+ [\gamma]_i $
{\ }and  {\ } $s_- : \NN^* \rightarrow \ZZ$ {\ } by {\ }
$(s_-)_j= \tilde u_j -\langle \gamma\rangle_j .$

Then $s_+,s_-$ are partitions and 
 
 $$H^q(G_r (V), S_w Q \otimes {\wedge}_u S) =
 \delta_{q, i(w, \tilde u)}{\mathcal S}_{\psi(w, \tilde u)}V, $$
with
\begin{description}

\item [(i)]  $\psi(w, \tilde u)  =  (s_+, s_-)^{\geq },$

\item [(ii)]   $i(w, \tilde u)  =  |s_+|-\sum_{i\in I(r)}w_i.$
\end{description}

\end{lem}

{\it Proof:}

Since the  cohomology group is given by Corollary 5.4,  we only need   
to investigate the combinatorics.

Let $a = (w, \tilde u )$ be admissible. With the above notations the
 set of elements in 
$(a-I(d))$ 
is $\{ \{\alpha_i\}_{i=1, \ldots, r},
\{ \beta_j\}_{j=1, 
\ldots,d-r} \}$, where $\alpha_i$ {\ }is in position 
{\ }$i${\ } in
{\ }$(a-I(d))$  and in position{\ } $i+[\gamma]_i$ {\ }in {\ }
$(a-I(d))^{\geq}$. Thus the term in this position
in {\ } 
$(a-I(d))^{\geq}+I(d)$
is $\alpha_i+ i+[\gamma]_i=(s_+)_i$.

Similarly, {\ } $\beta_j${\ } is in the position{\ } $r+j${\ } in{\ }
$a-I(d),$  
{\ }and in the position{\ } 
$r+j-\langle\gamma\rangle_j$
in  $(a-I(d))^{\geq}$, such that the term in this position 
in {\ } 
$(a-I(d))^{\geq}+I(d)$ is $\beta_j+r+j-\langle\gamma\rangle_j
=(s_-)_j$.

For admissible $a$ the sequences 
$s_+$ and $s_-$ are weakly decreasing. Indeed,  

$(s_+)_i - (s_+)_{i+1} = w_i-w_{i+1}+ {\rm card}\{j{\ }|{\ }
\alpha_i>\beta_j>\alpha_{i+1} \}$ and 

$$ {\rm card} \{j \in \NN^*{\ }| {\ } w_i-i>\beta_j >w_{i+1}-(i+1)\}
\leq w_i-w_{i+1},$$
since {\ } $\beta_j ${\ } is a strictly decreasing integral
sequence, and
similarly for {\ } $s_-.$ 

Since $s_-$ converges to $0,$ it is a partition. The sequence $s_+$
has a finite number of terms, which are greater than the limit
of $s_-$. Thus $s_+$ is a partition, too.

Finally from $(s_+)_i =w_i+ [\gamma]_i $ we get (ii).
$\hfill{\square}$

For $(w,\tilde u)$  admissible, we will apply this lemma 
in the case where 
$w/\chi(u)$ is a skew partition of length ~ {\ } $r.$ 

In this case we denote by 
$$
\begin{array}{lll}
\Sigma_{\pm} &=& \{ (i,j)\in {\mathcal D}(w/\chi(u)){\ }|{\
}\beta_{1-j}\gtrless\alpha_i\}.\\
{ }& = & \{ (i,j)\in {\mathcal D}(w/\chi(u)) {\ }|{\ } 
\tilde u_{1-j}-(r+1-j)
\gtrless  w_i-i \}
\end{array}
$$

Note that $(w,\tilde u)$ is admissible, iff {\ }{\ }
$ {\mathcal D}(w/\chi(u))  = \Sigma _+ \cup \Sigma_-.$
We denote   $ {\mathcal D}(w/\chi(u))$ by $\Sigma.$

\begin{Def}

We call $(w/\chi(u)$ admissible iff $\Sigma  = \Sigma _+ \cup
\Sigma_-,$ and we call  
 $\Sigma _+, \Sigma_-$   the splitting of ${\mathcal D}(w/\chi(u))$

\end{Def}

Recall that we put $\tilde u_j = + \infty$
 for $j\leq 0$, such that all $(i,j)\in\Sigma$ with $j>0$
  belong to
$\Sigma_+$.

\begin{lem}

For $(w,\tilde u)$  admissible and $w/\chi(u)$  a skew partition 
of length ~  $r, $ we have
\begin{description}

\item [(i)]  $s_{+}  =  [\Sigma _+]$

\item [(ii)] $s_{-}  =  \langle\chi ^*(\Sigma _-)\rangle$

\item [(iii)] $i(w,\tilde u)  =  |u|-{{\rm card}}(\Sigma_-).$

\end{description}
\end{lem}

{\it Proof} :
$$ 
\begin{array}{lll}
(s_+)_i & = & w_i + {\rm card}{\ } \{ j\leq 0{\ }| {\ }\tilde 
u_{1-j}-1+j-r > w_i-i \}\\
{} & = & w_i + {\rm card}{\ } \{ j\leq 0{\ }| {\ }\beta_{1-j} > \alpha_i \}.$$
\end{array}
$$

If $j\leq -u_{r+1-i}$ one has  $1-j > u_{r+1-i}$, thus{\ } 
$\tilde u_{1-j} < r+1-i.$ 
Moreover since ${\mathcal D}(\chi(u))\subseteq {\mathcal D}(w), $ we have  
$ -u_{r+1-i}\leq w_i.$ This implies  $j \leq w_i, $ thus $\beta_{1-j} < 
\alpha_i, $
and we can write
$$(s_+)_i = w_i + {\rm card}{\ } \{ j \leq 0 {\ } | {\ }
 -u_{r+1-i} < j\ ,\ \beta_{1-j} > \alpha_i \}.$$

If $0 \geq j > w_i,$ we have $ {\ } 1-j \leq u_{r+1-i}.$ Thus 
 $\tilde u_{1-j} \geq r+1-i $ and  $\beta_{1-j} > \alpha_i.$
 
This  yields for any  $w_i \in \ZZ$

$$(s_+)_i = {\rm card}{\ } \{ j\in \ZZ
 {\ }|{\ }-u_{r+1-i} < j \leq w_i,\beta_{1-j} > \alpha_i \}
 = [\Sigma _+]_i.$$ as required.

Similarly we have 

$$(s_-)_j = \tilde u_j- {\rm card} {\ }\{ i\in I(r)
{\ } | {\ } \beta_j > \alpha_{r+1-i}\} =
\tilde u_j - \langle \gamma\rangle_j.$$

a) If $u_i < j$ we have $\tilde u_j < i,$ and since 
${\mathcal D}(\chi(u))
\subseteq {\mathcal D}(w), $ we have  
$ -u_i\leq w_{r+1-i}.$ this implies $\beta_j < \alpha_{r+1-i},$ thus

$\langle \gamma\rangle_j =  {\rm card} {\ }
\{ i\in I(r) {\ } | {\ }-u_i < 1-j, {\ } \beta_j > \alpha_{r+1-i} \}.$

\bigskip

b) If $j \leq -u_{r+1-i}, $ we have $j < u_i$, thus $\tilde u_j \geq i$ 
and
$\alpha_{r+1-i} < \beta_j. $

Now
$$ 
\begin{array}{lll}
\Sigma_{-}& = &\{i\in I(r){\ } | {\ } -u_{r+1-i} < j\leq w_i,{\ }
 \alpha_i>\beta_{1-j} \}\\
\chi ^*(\Sigma_{-})& = &
\{i\in I(r){\ } | {\ } -u_i < 1-j \leq w_{r+1-i} , {\ }  \alpha_{r+1-i} 
>\beta_j \}\\

\langle \chi ^* \Sigma_{-}\rangle_j & = &
 
{\rm card} {\ }\{ i\in I(r){\ } | {\ }-u_i < 1-j \leq w_{r+1-i}, 
{\ } \alpha_{r+1-i} >\beta_{j} \}\\

{}& = &{\rm card} {\ }\{ i\in I(r){\ } | {\ }-u_i < 1-j ,
 {\ } \alpha_{r+1-i} >\beta_{j} \}.
\end{array}
$$
 The last equality is due to  b). Now

$$ 
\begin{array}{lll} 
\langle \chi ^* \Sigma_{-}\rangle_j + \langle \gamma\rangle_j & = &
{\rm card} {\ }\{ i\in I(r){\ } | {\ }-u_i < 1-j,
 {\ } \alpha_{r+1-i} >\beta_{j} \}\\
 {} &+& {\rm card} {\ }\{ i\in I(r){\ } | {\ }-u_i < 1-j
 {\ } \alpha_{r+1-i} <\beta_j  \}\\
 {}& = & {\rm card} {\ }\{ i\in I(r){\ } | {\ }-u_i < 1-j \}\\
{}& = & {\rm card} {\ }\{ i\in I(r){\ } | {\ }i \leq \tilde u_j \} 
= \tilde u_j  
 \end{array}
$$

as required

Finally  ${{\rm card}}(\Sigma_+)+{{\rm card}}(\Sigma_-)
=\sum_{i\in I(r)}
(w_i-\chi(u)_i)$. Together with ${{\rm card}}(\Sigma_+)=|s_+|$
and the previous lemma, this yields the desired formula for
$i(w,\tilde u)$.

$\hfill{\square}$

\begin{lem}
We have 

\begin{description} 

\item [(i)]  $\psi(w, \tilde u) = \left(\ [\Sigma_+],{\ } \langle\chi ^*(\Sigma _-)\rangle\ \right)^{\geq}$

\item [(ii)]  $\langle \Sigma_-\rangle_j \leq [\Sigma_+]_i 
\quad for \quad (i,j)
\in \Sigma_-$

\item [(iii)]   $[\Sigma_+]_i \leq \langle \Sigma_-\rangle_j 
\quad for \quad (i,j)
\in \Sigma_+$

\end{description}

\end{lem}
\smallskip 

{\it Proof:}

The preceding lemmas gives{\ } (i). We prove {\ } ii), the proof of {\ } iii) is
similar. 

Since $[\Sigma_+]_i$ is decreasing, it suffices to
consider for given $j$ 
  $$(i,j) \in \Sigma_-, {\ }(i+1, j) \not\in \Sigma_-.$$
  
Determine {\ }$ j'$ {\ } such that {\ }$(i,j')\in \Sigma_-,{\ }
 (i,j'+1)\notin \Sigma_-.$ 

Since {\ }$ (r+1-\tilde u_{1-j}, j)$ is the element of smallest
height in column $j$ of $\Sigma_-${\ }
and $(i,j)$ the element of greatest height,
we have 

$\langle \Sigma_-\rangle_j = \tilde u_{1-j}-r+i
\leq \tilde u_{1-j'}-r+i \leq w_i-j' = [\Sigma _+]_i.$ {\ }

$\hfill{\square}$

Since every row of ${\mathcal D}(w,\chi(u))$ contains exactly one
row of $\Sigma_+$ and every column of $u$ yields exactly one
column of $\Sigma_-$, the length of $\psi(w, \tilde u)$ is the
number of non-vanishing terms of $[{\mathcal D}(w/\chi(u))]$ plus $u_1$.

\smallskip

{\bf {Example}}

\smallskip

$w = (5,4,3,2,-1,-2)$ 

$u = (7,7,4,3,3,1), {\ } \chi (u) = (-1,-3,-3,-4,-7,-7),$  

$ \tilde u = (6,5,5,3,2,2,2)$

$\alpha = (4,2,0,-2,-6,-8), {\ } \beta = (-1,-3,-4,-7,-9,-10,-11)$ 

$[\gamma] = (0,0,0,0,3,4), {\ }\langle \gamma\rangle = (3,2,2,1,0,0,0)$

$s_+ = [\Sigma_+] =  (5,4,3,3,2,2), {\ }
 s_- = \langle \chi ^* \Sigma_{-}\rangle = (3,3,2,2,2,2).$

The whole shaded area is $w/\chi (u), $ the set of vertical strips is 
$\langle \chi ^* \Sigma_{-}\rangle$, the set of horizontal strips is $[\Sigma_+] $ 

\begin{figure}[htb]\label{fig2}
\begin{center}
\epsfig{file=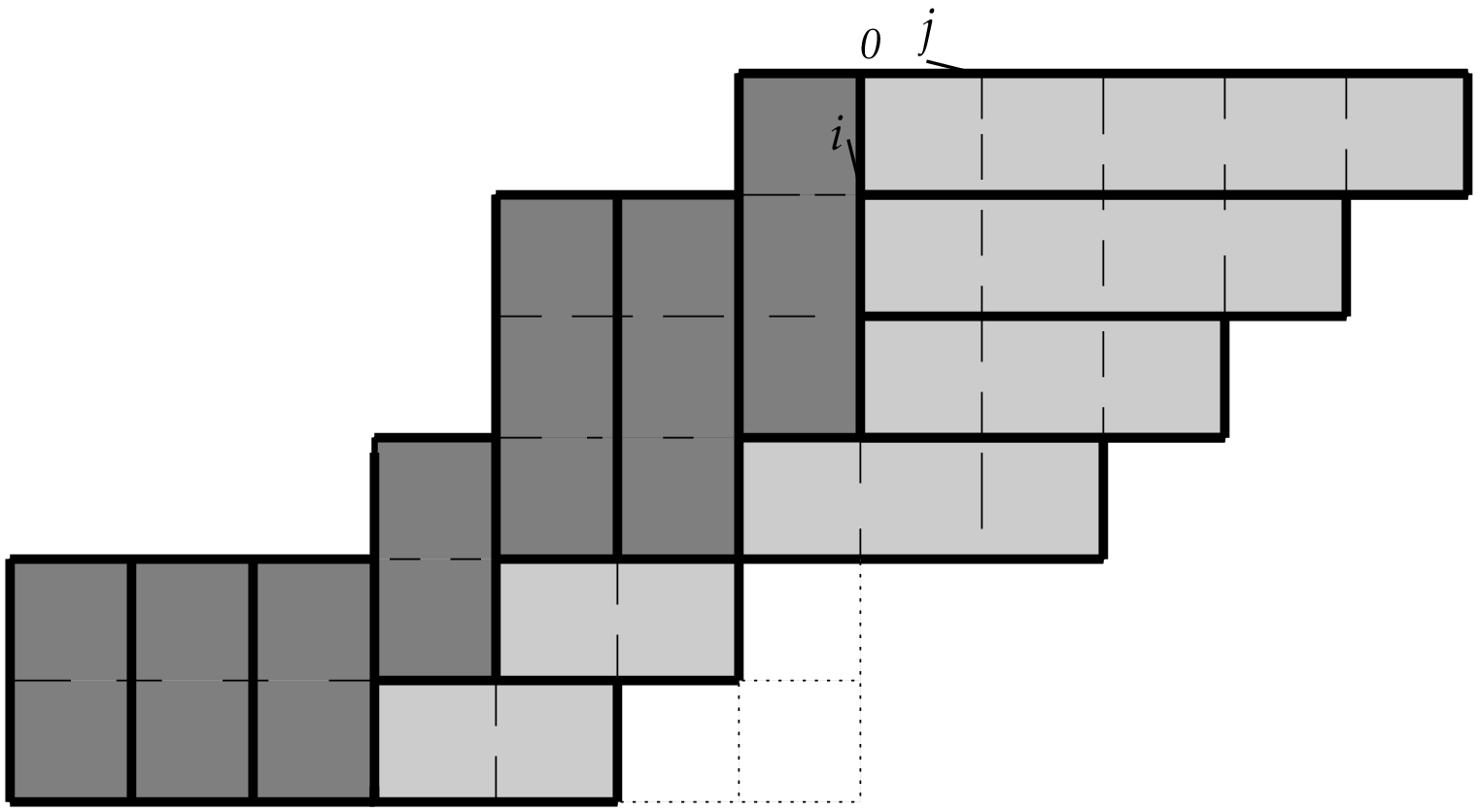,height=1.5in,width=3in,angle=0.0}
\caption{}
\end{center}
\end{figure}

\begin{lem}

Choosing the permutation of shortest length for the 
reordering of $( [\Sigma_+],{\ }
\langle\chi ^*(\Sigma _-)\rangle)$  yields a natural bijection 

$$\beta : \Sigma _+ \cup \Sigma_- \stackrel{\sim}\longrightarrow Y(
\psi(w,\tilde u )).$$

Then {\ } $\beta$ is height increasing.

\end{lem}

\smallskip 

{\it Proof:}

On {\ } $\Sigma_+$,{\ } the bijection $\beta$ is height increasing,
since {\ }$[\Sigma_+]$ {\ }  is weakly decreasing.
For  $(i,j) \in \Sigma_-$, the inequality derived above implies that
the reordering puts $\langle \Sigma_-\rangle_j$ after
$[\Sigma _+]_1,\ldots,[\Sigma _+]_i$. 

Thus {\ }  $\beta_1 (i,j) > i$.

$\hfill{\square}$

\begin{lem}

Let $u,v$ be partitions such that  $l(u)\leq r$ and $l(v) = r, $

$(w,b)\in LR(\chi(u), v)$, $(w,\tilde u)$ admissible, and 
$\rho=\psi(w,\tilde u)$. Then either

$1) {\ } u=0,{\ } \rho = v ,$  or

$2) {\ }|u|+i(w,\tilde u)+1 +||v||^2 \leq ||\rho||^2.$

\end{lem}

{\it Proof:}

We use induction on $u_1+r$, in other words on the length of
$\psi(w,\tilde u)$.  
For $u_1=0$, the statement is obviously true.

 Assume $u_1>0$. We use the splitting of ${\mathcal D}(w/ \chi(u))$
  into
$\Sigma_+,\Sigma_-$ introduced above.

Since $v_r>0$, we either have $\chi^\ast(1,u_1)\in \Sigma_-$
or $\chi^\ast(1,u_1)\in \Sigma_+$. The two cases will be treated
differently.

{\bf a)} For  $\chi^\ast(1,u_1)\in \Sigma_-$,
consider the partitions $u', v'$ given by

${\mathcal D}(u') = {\mathcal D}(u)
\cap (I(r)\times u_1)^c,$ and $ Y(v') = 
b^{-1}({\mathcal D}(w/ \chi (u'))),$  in other words if 
$\tilde u = ( \tilde u_1, \tilde u_2, \ldots, \tilde u_{u_1})$, 
with  $\tilde u_{u_1}=l, $ then 
$\tilde u' = ( \tilde u_1, \tilde u_2, \ldots, \tilde u_{u_1-1}).$
We denote by 
$x_j = (r+1-j, 1-u_1), j = 1,2,...l$ and $ L = \{x_j, j = 1,2,...l \}$

Let $b'$ be the restriction of $b$ to $Y(v')$.
By Remark 4.13 , $Y(v')$ is a Young diagram.
The Littlewood-Richardson rules
yield $(w,b')\in LR(\chi (u'),v').$
Note that $u'_1=u_1-1$ and $i(w,\tilde u)=i(w,\tilde u')$.
The skew partition $w/\chi (u')$ is admissible and
the splitting of ${\mathcal D}(w/ \chi (u'))$ is given by
$\Sigma'_+=\Sigma_+,{\ }$and   $\Sigma_- = \Sigma'_- \cup L$ 

By Lemma 4.5 and Remark 4.13, $Y(v)$ is obtained from $Y(v')$ by
 successive unions with the
preimages  of the set $L$, each union being
a Young diagram. Thus we have

$$||v||^2-||v'||^2=\sum_{j=1}^L (2c_{1}(x_j)-1),$$

Since $b$ is height increasing ie $c_{1}(x_j)\leq r+1-j$  this yields
$$||v||^2-||v'||^2\leq l(2r-l).$$
The length of $\rho=\psi(w,\tilde u)$ is $u_1+r$.
With $\rho'=\psi(w,\tilde u')$, lemma 5.9 yields 
$\rho=(\rho',l)$
and
$$||\rho||^2-||\rho'||^2 = l(2(u_1+r)-1).$$

By Lemma 5.8, we have $l =\langle\Sigma_-\rangle_{-u_1}\leq [\Sigma_+]_r.$

Since $[\Sigma_+]_r>0$ we have $v'_r>0$ and can use the
induction assumption. Thus

\begin{multline*}
|u|+i(w,\tilde u)+1+||v||^2-||\rho||^2\leq \\
l+ l(2r-l)- l(2(u_1+r)-1)+1 \\
=l(2-l-2u_1)+1\leq 0,
\end{multline*}
as required.

{\bf b)} For  $\chi^\ast(1,u_1)\in \Sigma_+$, the argument is similar.
We put

$${\mathcal D}(w'/\chi (u'))=
{\mathcal D}(w/\chi (u))\cap (\{r\}\times \ZZ)^c,{\ }{\mbox and} 
{\ }Y(v')= b^{-1}{\mathcal D}(w/\chi (u')),$$

in other words if $u = (u_1, u_2, u_3 \ldots )$ then 
$u' = (u_2, u_3 \ldots ).$

Let $b'$ be the restriction of $b$ to $Y(v')$.
The new partitions $v',w'$ have length $r'=r-1$.
The skew partition $w'/\chi (u')$ is admissible and
the splitting of ${\mathcal D}(w'/\chi (u'))$ is given by

$\Sigma'_-=\Sigma_-$.
Let $[\Sigma_+]_r=l$. We find
$$||v||^2-||v'||^2\leq l(2r-1),$$
$\rho=(\rho',l)$
and
$$||\rho||^2-||\rho'||^2 = l(2(u_1+r)-1).$$
Since $\langle\Sigma_-\rangle_{1-u_1}\neq 0$, we have
$u'_1 = u_1$ and  $v'_{r'} > 0$,
we  can use the induction assumption.
Altogether

\begin{multline*}
 |u|+i(w,\tilde u)+1+||v||^2-||\rho||^2\leq \\
2l + l(2r-1) -l(2(u_1+r)-1)\\
=2l(1-u_1)\leq 0.
\end{multline*}
$\hfill{\square}$

\smallskip

{\bf {Example of the situation of Lemma 5.10 case a)}}

For $u = (7,7,4,3,3,1)$  and  $ v = (9,8,6,5,5,3). $ 
We have $w/\chi (u)$ on the left hand side the same as the 
previous example with the same $[\Sigma_+],$ and 
$\langle \chi ^* \Sigma_{-}\rangle .$
The shaded area is $\tilde u_{u_1}.$ The partition 
$v'$ is the partition $v$ without the crossing boxes.

\begin{figure}[htb]\label{fig3}
\epsfig{file=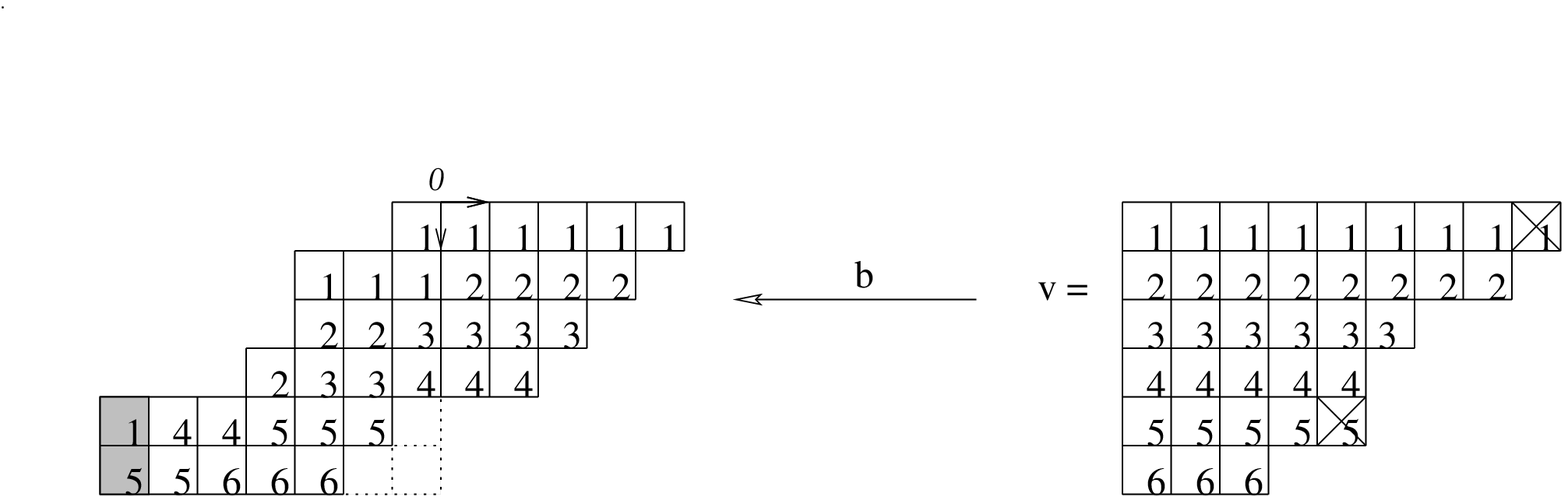,height=2in,width=5in,angle=0.0}
\caption{}
\end{figure}


\begin{lem}

Let $v$ be a partition of length $r$.
Then $$H^{p,q}(G_r(V), {\mathcal S}_v Q) =
 \bigoplus_{k\in K(p,q,r,v)} 
{\mathcal S}_{\rho(k)} V,$$

where $|\rho(k)| = |v|.$ We have 
\begin{description}

\item [(i)] $\forall k\in K(p,q,r,v)$, $\rho(k) \preceq v.$

\item [(ii)] Moreover, $\rho(k) = v$ only occurs for $p = q = 0,$
 where
 
  $$\mbox {for any partition }{\ }  v \qquad  H^0(G_r (V), 
  {\mathcal S}_v Q) =
   {\mathcal S}_v V.$$

\end{description}

\end{lem}

{\it Proof:}
     
It is well known that 
$$\Omega_{G_r(V)}^p = {\wedge} ^p(Q^*\otimes S) =\bigoplus_
{u\in{\sigma^p}}{\mathcal S}_u Q^*\otimes{\mathcal S}
 _{\tilde u}S,$$
where $\sigma^p$ is the set of partitions of weight $p$ and length $r$.                                                                
Then 

$$
\begin{array}{ll}
H^{p,q}(G_r(V),{\mathcal S}_v Q)&= \bigoplus _{u\in \sigma^p}
H^q(G_r (V), {\mathcal S}_{\chi (u)} Q\otimes
{\mathcal S}_v Q \otimes {\wedge}_u S)\\
&= \bigoplus _{u\in \sigma^p}  \bigoplus _{(w, b)\in LR(\chi (u), v)}
H^q(G_r (V), {\mathcal S}_w Q\otimes {\wedge}_u S)\\
&= \bigoplus _{u\in \sigma^p}  \bigoplus _{(w, b)\in LR(\chi (u), v)} 
\delta_{q, i(w,\tilde u)}{\mathcal S}_{\psi(w,\tilde u)}V.
\end{array}
$$
 
For each term on the right-hand side.
we have constructed a height increasing bijection 

$$\beta \circ b {\ } : {\ } Y(v)\stackrel{\sim} {\longrightarrow} 
Y(\psi (w,\tilde u)).$$
Thus $\psi (w,\tilde u) \preceq v.$
Since the length of $\psi (w,\tilde u)$ is at least $u_1$
plus the length of $v$ , it follows that
$\psi (w,\tilde u) = v $ implies $u = 0,$
 and thus $p = 0.$
 
$\hfill{\square}$

\begin{corol}
$H^{(p,q)}(G_r(V), det{\ } Q) = 0$ {\ } if {\ } $p \neq 0$ or $q \neq 0,$
\end{corol}

{\it Proof:}
   
there is no non-trivial partition strictly less 
than $v=(1,1,\ldots,1)$ of the same weight as $v.$ 
 
 This is also a result of Le Potier [LP1,(corol.1)]

$\hfill{\square}$

In the sequel, we will use the notation $R^{p,q}{\mathcal F}$ 
 for   $R^q\pi_*(\Omega ^p\otimes {\mathcal F})$.

\medskip

\subsection{ Proof of Theorem 5.1}$\phantom{0}$\\

\bigskip

For $p = 0$, Corollary 1 to Bott's theorem gives $i(a) = 0 = q$ 
and $\psi (a_s) = a_s.$ 

For  $p \neq 0$ we will use induction on $l$, the length of $s$.
Let us consider the Borel-le Potier (B-L) spectral 
sequence  associated
to $$\pi : Y = {\mathcal F}_s(V) \longrightarrow G_{s_{l}}(V) = X.$$
On $X$ we have the canonical quotient bundle $Q$ with fibres
$Q = V/V_{s_{l}}$. The fibres of $Y$ have the form
${\mathcal F}l_{s'}(Q_x)$, where $s'=(s_1,\ldots,s_{l-1})$.
On $Y$ we have canonical bundles $Q_i$ with fibres 
$V_{s_{i-1}}/V_{s_i}$.

As explained in the introduction,
the Leray spectral  sequence (L) associated to the projection  $\pi$, 
called $^{p',p}{\mathcal E}_L$ abuts to the ${\mathcal E}_1$ terms  
of the Borel-Le Potier (BL)  spectral sequence:

$$ {}^{p',p}{\mathcal E}^{q-j,j}_{2,L} \stackrel {L}
\Longrightarrow {} ^p{\mathcal E}^{p',q-p'}_{1, B}
\stackrel {BL}\Longrightarrow H^{p,q} ({\mathcal F}_s (V), Q^a).$$

We have 

$$^{p',p}{\mathcal E}^{q-j,j}_{2,L} = H^{p',q-j}(G_{s_l}(V),
\ R^{p-p',j}\ \pi_*(Q^a)).$$ 

On $Q$ we have flags 
$\{0\}\subset V_{s_{l-1}}/V_{s_l} \subset \dots V_{s_1}/ V_{s_l}\subset
V/V_{s_l}.$
For $ V_{s_{j}}/V_{s_l}= V'_{s_{j}}$ we have
$V'_{s_{j-1}}/V'_{s_{j}} = Q_j.$ We rewrite
$$
\begin{array}{lll}
Q^a &=& (\det \, Q_1)^{a_1}\otimes\ldots   \otimes (\det \, Q_l)^{a_l}\\ 
&= &(\det \,Q)^{a_l}\otimes (\det\, Q_1)^{a'_1}\otimes \ldots \otimes
(\det \, Q_{l-1})^{a'_{l-1}}
\end{array}
$$
where $a'_i = a_i-a_l$.

Setting $a'= (a'_1, \ldots, a'_{l-1}),$ we have

$$
\begin{array}{ll}
^{p',p}{\mathcal E}^{q-j,j}_{2,L} & = 
 H^{p',q-j}(G_{s_{l}}(V),(det Q)^{a_l}
\otimes R^{p-p',j}\ \pi_*(Q^{a'}))\\
& =  H^{p',q-j}(G_{s_{l}}(V),(det Q)^{a_l} \otimes H^{p-p',j}
({\it F}_{s_1,s_2 \dots s_{l-1}}(V),Q^{a'})). 
\end{array}
$$

For $p-p'=0$, the desired result follows from Corollary 5.4 
with $v=0$ and Lemmas 5.10 and 5.11. For $p'<p$ we have by
 the induction
assumption for the graded bundle $(Gr)$ associated 
to the higher direct
image, 

$$Gr\ R^{p-p',j}\ \pi_*(Q^{a'}) =
\bigoplus_{k\in K}{\mathcal S}_{\rho'(k)}\ Q,$$
where $\rho'(k)\preceq a'_{s'}$ for all $k\in K$.
Consequently, 

$$(\det  Q)^{a_l} \otimes {\mathcal S}_{\rho'(k)}Q
= {\mathcal S}_{\rho''(k)}Q \qquad  {\mbox  with} \qquad  \rho''(k)\preceq a_s.$$

Since $|\rho'(k)|=|a'_{s'}|$, we have $|\rho''(k)|=|a_s|$.

The bundle $R^{p-p',j}\ \pi_*(Q^{a'})$ on $X$ is a homogeneous
$Gl(V)$-bundle, thus specified by a representation of the
stabilizer of a point in $X$. Since $G(V_{s_l})$ acts trivially 
on the fibres, its representation factorizes through that
of $Gl(Q)$. By the Schur lemma [FH],
such a representation is reducible and is given by a Schur functor,
which implies 

$$R^{p-p',j}\ \pi_*(Q^{a'}) \simeq Gr\ R^{p-p',j}\ 
\pi_*(Q^{a'}).$$

Finally, again by the induction assumption,

$$p-p'+j+1+||a'_{s'}||^2\leq ||\rho'(k)||^2$$
for all $k$ in the  set of subscripts $K$. By Lemma 5.10 we have
$$p'+q-j+||\rho'(k)||^2 +|a_l|^2\leq ||\rho(i)||^2$$
for all $i\in I(p,q,s,a)$. Since $||a_s||^2 = 
||a'_{s'}||^2 + |a_l|^2$,
the result follows from Lemmas 5.10 and 5.11.
$\hfill{\square}$
\bigskip

\section{ Proof of Theorems 1.1 and Corollary 1.2}

\subsection{ Proof of Theorem 1.1} $\phantom{0}$\\

For a partition   
$R =(r_1, \ldots, r_m)$, we take 
{\ } {\ }$a= \tilde R^>, {\ }s = (s_1, \ldots , s_l) = R^<\ ,$
such that $a_s = \tilde R.$

Let $Y  = \mathcal F_{s_1,\ldots ,s_l}(E) $ and $\mathcal L = Q^a$
We use induction on $R$ with respect to the dominance partial
order.

In the setup discussed in the introduction, we consider
the Borel- Le Potier spectral sequence given by
$X,Y,{\mathcal L},P$, where ${\mathcal L}=Q^a$. Its ${\mathcal E}_1$
terms can be evaluated by a Leray spectral sequence, for which 
$$^{p,P}{\mathcal E}_{2,L}^{q-j,j} =
H^{p,q-j} (X, R^{P-p,j}\pi_*(Q^a)), $$

For $P-p = 0, $  we have {\ }
$R^{0,j}\pi_*(Q^a) = \delta _{j,0}{\ } {\wedge} _rE$ {\ }
by Theorem 5.1.
This implies that for {\ } $P-p = 0$ {\ } the Leray spectral sequence
degenerates at ${\mathcal E}_2$,
and {\ }  $^P{\mathcal E}_{1,B}^{P,q-P} = H^{P,q}(X, {\wedge}_r E).$

We have 

$$
\xymatrix{ ^P{\mathcal
E}_{1,B}^{p-p_1,(q-1)-(p-p_1)}\ar[r]^{\phantom{xxxxxxx}
d_{1,B}}&\ldots \ar[r]^{\hspace*{-2.5em}d_{1,B}}& ^P{\mathcal
E}_{1,B}^{p-1, 
(q-1)-(p-1)} \ar[r]^{\phantom{xxxx}d_{1,B}}& {}^P{\mathcal E}_{1,B}^{p,q-p}
\ar[r]^{\phantom{xx}{d_{1,B}}} & 0\\ 
\ldots \ar[urrr]_{d_{p_1,B}}&&&&}
$$

In order to show that $^P{\mathcal E}_{1, B}^{P,q-P}$ is a subfactor
of the limit group $H^{P,q}(Y,Q^a),$
we must prove that the sources of  {\ }  $\stackrel
{d_{p_1,B}}\longrightarrow  { }^P{\mathcal E}_{1,BL}^{p,q-p}$
vanish for each $p_1\not= 0$. 
 
Now each  group{\ } $^P {\mathcal E}_{1,B}^{p-p_1, (q-1)-(p-p_1)}$ 
 is a subquotient of

 $$H^{P-p_1, q-q_1-1}(X, R^{p_1,q_1}(Q^a)),$$
 
where by Theorem 5.1,  
$$\rm Gr\ R^{p_1,q_1}\pi_*(Q^a) = \bigoplus{\mathcal S}
_{\rho(i)}(E)$$
with ${\rho(i)} \prec\tilde R$, since $p_1\neq 0$.
By Lemmas 3.5 and 3.6 the vector bundles ${\mathcal S}_{\rho(i)}(E)$ 
are ample.
Moreover, $p_1+q_1+1 \leq |\rho|^2-|R|^2$ by Theorem 5.1.
Together with the assumption

$$ P+q-n> \sum_{i=1}^m r_i(d-r_i){\ } {\ } {\ } {\ }(*)$$
and $|\rho|=|R|$ this yields
$P-p_1+q-q_1-1-n > \sum_{i=1}^m\rho_i(d-\rho_i)$. 
Since ${\rho(i)} \prec\tilde R$ we can use the induction assumption,
such that indeed
$H^{P-p_1, q-q_1-1}(X, Gr \ R^{p_1,q_1}(Q^a))$ and consequently

$H^{P-p_1, q-q_1-1}(X, R^{p_1,q_1}(Q^a))$
and $^P {\mathcal E}_{1,B}^{p-p_1, (q-1)-(p-p_1)}$ are zero.

It may well be that $R^{p_1,q_1}(Q^a)$ is isomorphic to the
associated graded bundle $Gr \ R^{p_1,q_1}(Q^a)$, as in 
the Grassmannian case discussed above. We did not investigate this 
question, since it is not necessary for our proof.

We have shown that under the condition (*) $^P{\mathcal E}_{1,B}^{P,q-P} = 
H^{P,q}(X, {\wedge}_r E)$ is a subquotient of $H^{P,q}(Y,Q^a)$.

Now $Q^a$ is ample by 

 \smallskip

\begin{lem}

Let $a=({a_1,a_2,\ldots})$ a strictly decreasing partition, and {\ } $s$
as above.

If {\ }${\mathcal S}_{a_s}E$ {\ } is ample over $X$ , then $Q^a$ is ample
over $Y$.

\end{lem}

{\it Proof:}  See [D1], Lemmas 2.11 and 4.1.  $\hfill{\square}$

\smallskip

To conclude the proof of Theorem 1.1, the group $H^{P,q}(Y,Q^a)$ 
vanishes under 
the condition of Kodaira-Akizuki-Nakano theorem but the condition 
(*) implies this condition.

\subsection{ Proof of Corollary 1.2} $\phantom{0}$\\

We need the following 
\begin{lem}
$\otimes_{i=1}^l S^{k_i}E\otimes _{j=1}^m 
\wedge^{s_j} E$ is ample $\Longleftrightarrow \mathcal S_{\lambda} E $ is ample,
where $\tilde \lambda = (s_1,s_2,\ldots , s_m,  \underbrace{1,\ldots, 1}
_{k_1\  {\rm times}}, \underbrace{1,\ldots, 1}
_{k_2\  {\rm times}},\ldots,\underbrace{1,\ldots, 1}
_{k_m\  {\rm times}}).$  
\end{lem} 

{\it Proof:}    
  
$\mathcal S_{\lambda} E $ is direct summand of  $\otimes_{i=1}^l
 S^{k_i}E\otimes _{j=1}^m 
\wedge^{s_j} E.$ 
  
Conversely for any $\mathcal S_{\lambda_i}E$ direct summand  of   $\otimes_{i=1}^l
 S^{k_i}E\otimes _{j=1}^m 
\wedge^{s_j} E,$ we have $ \lambda  {\ }\succeq {\ }\lambda  _i,$ then we use Lemma 3.5.   
$\hfill{\square}$
  
Now

 Among $\mathcal S _{\lambda_i} E$ where $ {\lambda_i} = (r_1,r_2, \ldots )$, 
that appear as direct summands of $\otimes_{i=1}^l S^{k_i}E\otimes _{j=1}^m 
\wedge^{s_j} E, $ the optimal condition $\sum\limits_{i=1}^m r_i(d-r_i)$ is
obtained by  $\lambda.$
 
This concludes the proof of Corollary 1.2 
$\hfill{\square}$

\end{document}